\documentclass[11pt]{article}

\usepackage{tgpagella}
\usepackage[utf8]{inputenc}
\usepackage[T1]{fontenc}
\usepackage{amssymb}
\usepackage{amsfonts}
\usepackage{amsmath}
\usepackage{amsthm}
\usepackage{mathrsfs}
\usepackage{comment}
\usepackage[normalem]{ulem}
\usepackage{hyperref}

\newtheorem{theorem}{Theorem}

\newtheorem{lemma}[theorem]{Lemma}
\newtheorem{fact}[theorem]{Fact}
\newtheorem{cor}[theorem]{Corollary}

\theoremstyle{definition}
\newtheorem{definition}[theorem]{Definition}
\newtheorem{remark}[theorem]{Remark}

\newtheorem{claim}{Claim}

\newcommand{\FS}{\textnormal{FS}}
\newcommand{\PA}{\textnormal{PA}}
\newcommand{\KF}{\textnormal{KF}}
\newcommand{\VF}{\textnormal{VF}}

\newcommand{\CD}{\textnormal{CD}}

\newcommand{\CT}{\textnormal{CT}}
\newcommand{\df}[1]{\textbf{#1}}
\newcommand{\num}[1]{\underline{#1}}

\newcommand{\ElDiag}{\textnormal{ElDiag}}
\newcommand{\val}[1]{{#1}^{\circ}}
\newcommand{\Val}{\textnormal{Val}}

\newcommand{\dom}{\textnormal{dom}}

\newcommand{\qcr}[1]{\ulcorner #1 \urcorner}

\newcommand{\LPA}{\mathscr{L}_{\PA}}
\newcommand{\form}{\textnormal{Form}}
\newcommand{\Term}{\textnormal{Term}}
\newcommand{\Sent}{\textnormal{Sent}}
\newcommand{\Var}{\textnormal{Var}}
\newcommand{\ClTerm}{\textnormal{ClTerm}}

\newcommand{\ClTermSeq}{\textnormal{ClTermSeq}}

\newcommand{\Comp}{\textnormal{Comp}}

\newcommand{\Asn}{\textnormal{Asn}}

\newcommand{\synt}[1]{\mathsf{#1}}

\title{Classical determinate truth without induction}
\author{Bartosz Wcisło}

\begin{document}
\maketitle

\begin{abstract}
	Fujimoto and Halbach had introduced a novel theory of type-free truth $\CD$ which satisfies full classical compositional clauses for connectives and quantifiers. Answering their question, we show that the induction-free variant of that theory is  conservative over Peano Arithmetic.
\end{abstract}

\section{Introduction}

Formal truth theory studies the concept of truth and various philosophical questions surrounding that notion using the methods of formal logic. An important part of that subject studies such questions by investigating \emph{axiomatic truth theories} obtained by adjoining to a sufficiently strong base theory (typically Peano Arithmetic, $\PA$), a fresh predicate $T(x)$ with the axioms stipulating that it satisfies certain truth-like properties, for instance, Tarski's compositional clauses. See \cite{halbach} for a standard introduction to the subject. 

One of the main topics of this study is the possibility of formulating well-defined and robust principles governing \emph{self-referential} truth predicates, i.e., predicates in which the constructed truth predicate can be applied not only to the sentences from the language of the base theory or some syntactically restricted fragment thereof, but -- ideally -- to the whole language. Such theories are also referred to as \emph{untyped} or theories of \emph{untyped truth}. Probably the most prominent example is the theory $\KF$ initially introduced introduced in \cite{feferman}, axiomatising Kripke's theory of truth, as laid out in \cite{kripke}. Other prominent examples are $\FS$, isolated in \cite{friedman_sheard} and $\VF$ (see \cite{cant}), axiomatising the supervaluational concept of self-referential truth.

In \cite{fujimoto_halbach_cd1}, another untyped theory of truth called $\CD$ has been introduced. Its main feature is that fully classical compositional axioms apply to all the sentences in its language with the only restrictions applying to atomic formulae of the form $T(s)$, where in order to apply the usual Tarskian clause we additionally require that the formula in question be determinate which, in turn, is introduced as a separate primitive predicate.

The theory aims to capture a very natural intuition: there are some sentences whose truth value is \emph{determined}. Even though they may  contain the truth predicate, we can uncontroversially ascribe to them a well defined truth value. An example of such a sentence is: ``It is true that it is true that $0=0$.'' Even though the truth predicate nontrivialy appears in this sentence, we generally would not find it problematic to say that it is simply true. On the other hand, the liar's sentence seems to really lack any well-defined truth value. In $\CD$, this distinction gets promoted to a primitive predicate: we stipulate that there is a class of determinate sentences, satisfying certain closure conditions, and that the truth predicate fully commutes with all the connectives and quantifiers, and admits full disquotation for the atomic sentences, including sentences of the form $T(\phi)$ provided that they are determinate. 

The theory $\CD$ turns out to be proof-theoretically equivalent to Kripke--Feferman theory $\KF$. The authors conjectured that the theory $\CD^-$ in which induction is restricted only to formulae in the base language is in fact conservative.

In our note, we confirm this conjecture. We will show that in every countable recursively saturated model of $\PA$, one can find an expansion to a model of $\CD^-$. Our proof closely resembles the original argument for the consistency of $\CD$, where in some crucial places we use techniques typical for the study of nonstandard models of arithmetic and truth predicates.

\section{Preliminaries} \label{sec_prelim}

All the theories of truth which we will be considering are formulated over Peano Arithmetic, $\PA$. Crucially, in $\PA$ we can formalise syntactic notions thus we can speak of terms, formulae, or proofs. The details are covered in a number of sources, for instance \cite{kaye} or \cite{hajekpudlak}. In order to make the text smoother, we will only occasionally comment on the formulae expressing the formalised syntactic concepts and instead present them all in Appendix \ref{sec_app_syntax}.

Throughout the paper, we will adopt a number of conventions to facilitate the reading:

\begin{itemize}
	\item We will conflate the formulae denoting syntactic operations and write them like the operations themselves. For instance, $P(\phi \wedge \psi)$ should be eliminated as saying: ``there exists a formula $\eta$ such that $\eta$ is the unique conjunction of $\phi$ and $\psi$, and $P(\eta)$ holds.''
	\item More generally, if a formula defines a syntactic operation, we will sometimes treat it as if it were a new term. For instance, if $\Val(x,y)$ is a formula expressing that the formally computed value of an arithmetical term $x$ is $y$, we will instead denote this value by $\val{x}$ and use it as if it were an independent expression.
	\item We will treat the formulae denoting syntactic notions as if they were denoting sets. For instance, we will be mostly writing $\phi \in \Sent$ instead of $\Sent(\phi)$.
	
\end{itemize}

Let us now introduce the theory which is the central subject of this article.

\begin{definition}
	By $\CD^-$ we mean the theory in arithmetical language with two additional predicates $T$ and $D$ extending $\PA$ with the following axioms:
	\begin{itemize}
		\item[$\synt{T1}$] $\forall s,t \in \ClTerm_{\LPA} \ T(s=t) \equiv \val{s} = \val{t}$. 
		\item[$\synt{T2}$] $\forall x \ Dx \rightarrow T D \num{x}$
		\item[$\synt{T3}$] $\forall x \ Dx \rightarrow \left(Tx \equiv TT\num{x}\right).$
		\item[$\synt{T4}$] $\forall \phi \in \Sent \ T\neg \phi \equiv \neg T \phi$. 
		\item[$\synt{T5}$] $\forall \phi \in \Sent \ T(\phi \wedge \psi) \equiv \left(T \phi \wedge T \psi\right).$
		\item[$\synt{T6}$] $\forall \phi, v \Big(\forall v \phi \in \Sent \rightarrow T \forall v \phi \equiv \forall x T\phi(\num{x})\Big).$ 
		\item[$\synt{D1}$] $\forall s, t \in \ClTerm_{\LPA} \ D(s=t).$
		\item[$\synt{D2}$] $\forall x \ DT(\num x ) \equiv D(x)$. 
		\item[$\synt{D3}$] $\forall x \ DD(\num{x}) \equiv D(x)$.
		\item[$\synt{D4}$] $\forall \phi \in \Sent \ D \neg x \equiv D x. $
		\item[$\synt{D5}$] $\forall \phi, \psi \in \Sent \ D(\phi \wedge \psi) \equiv \left((D(\phi) \wedge D(\psi)) \vee (D(\phi) \wedge T \neg \phi) \vee (D(\psi) \wedge T \neg \psi) \right).$  
		\item[$\synt{D6}$] $\forall \phi, v \ \Big((\forall v \phi) \in \Sent  \rightarrow D(\forall v \phi) \equiv \left( \forall x D \neg \phi(\num{x}) \vee \exists y D \phi(\num{y}) \wedge T \neg \phi(\num{y}) \right) \Big).$
		\item[$\synt{R1}$] $\forall \bar{s}, \bar{t} \in \ClTermSeq_{\LPA} \forall \phi \ \Big(\bar{\val{s}} = \bar{\val{t}} \rightarrow T\phi(\bar{s}) \equiv T \phi(\bar{t}) \Big).$
		\item[$\synt{R2}$] $\forall \bar{s}, \bar{t} \in \ClTermSeq_{\LPA} \forall \phi \ \Big(\bar{\val{s}} = \bar{\val{t}} \rightarrow D\phi(\bar{s}) \equiv D \phi(\bar{t}) \Big).$
	\end{itemize}
\end{definition}

Above, $\Sent,$ refers to (the formal definition of the codes of) the sentences in the full language, including $T$ and $D$ predicates. By $\val{s}$, we mean the formally computed value of the term $s$, when we write $\bar{\val{s}}$, we mean the sequence of values of the terms of the sequence of terms $\bar{s}$. If $x$ is an arbitrary element, $\num{x}$ is some canonically chosen numeral with the value $x$. The formulation of certain axioms is slightly different than the one in \cite{fujimoto_halbach_cd1} in that we do not quantify over arbitrary terms, but only over numerals. However, this is not really an issue, since we adopt regularity axioms anyway. These are the clauses $\synt{R1}$ and $\synt{R2}$ which say that the determinatess or truth of a sentence does not change if we replace its terms with any other terms of the same value. This assumption comes handy in the study of induction-free fragments of truth predicates and is clearly harmless when our goal is to prove conservativity.

If we add to $\CD^-$ the full induction scheme for the language containing the predicates $D,T$, we obtain a theory called $\CD$ which was the focus of study of the original article \cite{fujimoto_halbach_cd1}. The crucial results of their paper can be summarised as follows:

\begin{theorem}[Fujimoto--Halbach] \label{th_fujimoto_halbach_cd}
	The theory $\CD$ is consistent. Moreover, its arithmetical part is identical with the arithmetical part of Kripke--Feferman theory $\KF$.
\end{theorem}
The theory $\KF$ mentioned above is an extremely important theory of self-referential truth, which we have mentioned in the introduction.


\subsection{Models of arithmetic and truth theories}

In our paper, we will make the substantial use of nonstandard models of arithmetic. We assume that the reader is familiar to some extent with the topic. All the relevant notions and facts (with full proofs) can be found in \cite{kaye}.

As we have already mentioned in the introduction, our conservativity proof closely follows the original consistency proof of Fujimoto and Halbach. In the original paper, they have obtained a model of $\CD$ via an inductive construction over $\mathbb{N}$. In our paper, we will try to repeat their construction over an arbitrary model of $\PA$. On of the obvious obstacles is that in some steps of their proof, we are considering the set of true sentences over a given structure $(\mathbb{N},A)$. This step never yields a predicate satisfying Tarskian compositional clauses, when applied to a nonstandard model, since in this manner we fail to assign any truth value to nonstandard sentences $\phi$, thus violating the compositional axiom for the negation.

One fact that we will crucially use throughout the whole paper is resplendence of the countable recursively saturated models of $\PA$. 

\begin{fact} \label{fact_resplendence}
	Let $M \models \PA$ be countable and suppose than an expanded structure $(M,A)$ is recursively saturated. Suppose that $U$ is a recursive theory in an extended language consistent with $\ElDiag(M,A)$. Then there exist relations $R_{i}, i<\omega$ on $M$ such that the expansion $(M,A,R_i)_{i <\omega}$ satisfies the theory $U$.
\end{fact}
The above fact (in a somewhat weaker version) is called the \df{resplendence} of countable, recursively saturated models. 

In \cite{enayatvisser2}, the authors have introduced a new elegant method for proving conservativeness of the the classical stratified Tarskian truth predicate satisfying pure compositional clauses which since then has become an indispensable tool for the field. Let us phrase an immediate corollary of their results which we will make a crucial use of in our paper.

\begin{definition}
	Let $\CT(D,T,T')$ be a conjunction of the following clauses, expressing that a predicate $T'$ satisfies compositional conditions for the language with the predicates $D,T$:
\begin{itemize}
	\item $\forall s,t \in \ClTerm_{\LPA} \ \Big(T'(s=t) \equiv \val{s} = \val{t}\Big).$
	\item $\forall x \Big(T'D(\num{x}) \equiv D(x)\Big).$
	\item $\forall x \Big(T'T(\num{x}) \equiv T(x)\Big).$
	\item $\forall \phi \in \Sent \ T'\neg \phi \equiv \neg T' \phi$. 
	\item $\forall \phi \in \Sent \ T'(\phi \wedge \psi) \equiv \left(T' \phi \wedge T' \psi\right).$
	\item $\forall \phi, v \Big(\forall v \phi \in \Sent \rightarrow T' \forall v \phi \equiv \forall x T'\phi(\num{x})\Big).$ 
	\item $\forall \bar{s}, \bar{t} \in \ClTermSeq \forall \phi \ \Big(\bar{\val{s}} = \bar{\val{t}} \rightarrow T'\phi(\bar{s}) \equiv T' \phi(\bar{\val{t}}) \Big).$
\end{itemize}
\end{definition}

\begin{fact} \label{fct_ctminus_w_wiekszych_jezykach}
	Suppose that $(M,D,T)$ is an arbitrary model and $M \models \PA$. Then there exists $(M^*,D^*,T^*) \succeq (M,D,T)$ and $T' \subset M'$ such that $(M^*,D^*,T^*,T') \models \CT^-(D^*,T^*,T')$.
\end{fact}
As we have mentioned, the Fact follows by known arguments for conservativity of $\CT^-$, for instance those of Enayat and Visser, simply by proving it for larger languages. We are not aware of a formulation of that statement in the literature that exactly matches our situation. However, the proof of a related and slightly more complicated result presented in Appendix A, can be easily adapted to obtain the above Fact, as per Remark \ref{rem_enayat_visser_bigger_languages}.

\begin{cor} \label{cor_ctminus_in_respendent_models}
	Let $M \models \PA$ and assume that $(M,D,T)$ is a countable recursively saturated model. Then there exists $T' \subset M$ such that $(M,D,T,T') \models \CT^-(D,T,T')$.
\end{cor}
\begin{proof}
	This is an immediate corollary of Facts \ref{fact_resplendence} and \ref{fct_ctminus_w_wiekszych_jezykach}.
\end{proof}

\section{The main result}

The rest of the paper is devoted to the proof of the main theorem.

\begin{theorem} \label{th_cdminus_conservative}
	$\CD^-$ is conservative over $\PA$. 
\end{theorem}

\begin{proof}
	Fix a countable recursively saturated model $M \models \PA$. We will find $T,D \subset M$ such that $(M,T,D) \models \CD^-$. 
	
	First, let us introduce some notation. For a pair of predicates $D,T$, let $\mathscr{D}(D,T)(x)$ be the disjunction of the following clauses:
	\begin{itemize}
		\item $\exists s,t \in \ClTerm_{\LPA} \ x = (s=t)$.
		\item $\exists t \in \ClTerm_{\LPA} \ x = Tt \wedge D (\val{t})$.
		\item $\exists t \in \ClTerm_{\LPA} \ x = Dt \wedge D(\val{t}).$
		\item $\exists \phi \in \Sent \ x= \neg \phi \wedge D(\phi).$ 
		\item $\exists \phi, \psi \in \Sent \ x = (\phi \wedge \psi) \wedge \big((D(\phi) \wedge D(\psi)) \vee (D(\neg \phi) \wedge T(\neg \phi) ) \vee (D(\neg \psi) \wedge T(\neg \psi)) \big).$
		\item $\Sent(x) \wedge \exists v, \phi \ (x= \forall v \phi) \wedge \big((\forall y D \phi(\num{y})) \vee (\exists y D \neg\phi(\num{y}) \wedge T \neg \phi(\num{y})) \big).$
	\end{itemize}

	Notice that the above definition is slightly modified in comparison to \cite{fujimoto_halbach_cd1}: in the clauses for conjunction and the universal quantifier, we explicitly require the negation of the conjuncts (or of an instance of substitution) to be determined. This allows the bookkeeping in certain subsequent arguments to follow in a smoother manner.

	Now, consider the theory $\Theta$ which formalises Fujimoto--Halbach construction of a model of $\CD$. The language of the theory $\Theta$ consists of:
	\begin{itemize}
		\item The arithmetical language.
		\item A binary predicate $T(i,x)$ which we will denote $T_i(x)$ and think of as a family of unary predicates.
		\item A binary predicate $D(i,x)$ which we will denote $D_i(x)$ and think of as a family of unary predicates.
		\item A unary predicate $I(x)$ which will denote the domain of the indices $i$. 
	\end{itemize}
	
	We endow $\Theta$ with the following axioms:
	\begin{itemize}
		\item $\PA$.
		\item $\forall x \neg D_0(x)$.
		\item $\forall x \neg T_0(x)$.
		\item $\forall i \ \Big[ i+1 \in I \rightarrow \forall x \Big(D_{i+1}(x) \equiv \mathscr{D}(D_i,T_i)(x)\Big)\Big].$
		\item $\forall i \ \Big[i+1 \in I \rightarrow \CT(D_i, T_i,T_{i+1}\Big)\Big].$
		\item $\forall x,y \Big(x<y \wedge I(y) \rightarrow I(x)\Big)$.
		\item $I(n), n \in \omega$. 
	\end{itemize}
	
	\begin{claim}
		 The theory $\Theta$ is consistent. 
	\end{claim}

By compactness, it is enough to check whether it is consistent if we interpret $I$ as a (standard) finite interval $\{0,1, \ldots, n\}$. 
	
	Fix any $n \in \omega$. Take any countable recursively saturated model $M \models \PA$ and iteratively define the predicates $D_i,T_i: i < n$. Start from $D_0 = T_0 = \emptyset$. For $i$ such that $i+1 < n$, let:
	\begin{itemize}
		\item $D_{i+1}$ be defined from $D_i, T_i$ via the predicate $\mathscr{D}$. 
		\item $T_{i+1}$ be an arbitrary subset of $M$ such that $(M,D_i,T_i, T_{i+1})$ is a recursively saturated model satisfying $\CT(D_i,T_i,T_{i+1})$.  
	\end{itemize}
	We maintain throughout the construction that all the models obtained are recursively saturated. In the step for $T_{i+1}$ we use Corollary \ref{cor_ctminus_in_respendent_models}, resplendence of countable recursively saturated models and the fact that the recursive saturation itself can be expressed with a first-order theory. This concludes the proof of the claim.
	
	Now, fix an arbitrary countable recursively saturated model
	\begin{displaymath}
		(M,I,T^*,D^*) \models \Theta
	\end{displaymath}
	and let
	\begin{eqnarray*}
		D_{\omega} & := & \bigcup_{i \in \omega} D_i \\
		T_{\omega} & := & \bigcup_{i \in \omega} (T_i \cap D_i).
	\end{eqnarray*} 

Note that we are taking unions only over the indices $i$ in the standard part of the model. We will show that the model $(M,D_{\omega},T_{\omega})$ satisfies a number of desirable properties and then proceed to the construction of the actual model of $\CD^-$. 
		
\begin{claim} \label{cl_monotonicity}
	For arbitrary $k \leq n< \omega$, the following relations hold:
	\begin{itemize}
		\item $D_k \subseteq D_{n}$. 
		\item $D_k \cap T_k = D_{k} \cap T_{n}$. 
	\end{itemize}
\end{claim}		

\begin{proof}[Proof of Claim \ref{cl_monotonicity}]
	We show both claims by simultaneous induction on $(n,k)$ ordered lexicographically. Both claims are obvious for $k=0$ and an arbitrary $n$. We will show that they hold for $(n+1,k+1)$, assuming that they hold for all smaller pairs.
	
	First, fix any $\phi \in D_{k+1}$. We will prove by case distinction depending on the syntactic shape of $\phi$ that $\phi \in D_{n+1}$. Most of the cases could be simply summarised by saying that $\mathscr{D}$ is a positive operator. However, there is a small wrinkle about the case for conjunction and the universal quantifier, so we describe it in more detail.
	
	\paragraph*{Atomic case -- arithmetical formulae} Suppose that there exist $t,s \in \ClTerm_{\LPA}$ such that $\phi = (s=t)$. Then $\phi \in D_{n+1}$ by definition. 
	
	\paragraph*{Atomic case -- determinacy} Suppose that there exists $\psi \in D_k$ such that $\phi = D \psi$. Then by induction hypothesis, $\psi \in D_n$, and so $\phi \in D_{n+1}$. 
	
	\paragraph*{Atomic case -- truth} Suppose that $\phi = T\psi$ for some $\psi \in D_k$. Again, this implies by induction hypothesis that $\psi \in D_n$ and $\phi \in D_{n+1}$.
	
	\paragraph*{Negation} Suppose that $\phi = \neg \psi$ for some $\psi \in D_k$. Like in the previous, analogous cases, $\psi \in D_n$, and so $\phi \in D_{n+1}$.
	
	\paragraph*{Conjunction} This case is slightly different from the previous ones. Suppose that $\phi = \psi \wedge \eta$. 
	
	Then we have two cases. If $\psi, \eta \in D_k$, then we can conclude that $\psi \wedge \eta \in D_{n+1}$, like in the previous cases. So suppose that one of the formulae is determinate and false. Say, $\neg \psi \in D_k$ and $\neg \psi \in T_k$. Then by induction hypotheses, $\neg \psi \in D_k \cap T_k = D_k \cap T_n \subseteq D_n \cap T_n$, so $\psi \wedge \eta \in D_{n+1}$. 
	
	\paragraph*{Universal quantifier} This case is similar to the previous one. Suppose that $\forall v \phi(v) \in D_{k+1}$. If for all $x \in M$, $\phi (\num{x}) \in D_k$, then by induction hypothesis, they are also in $D_n$, so $\forall v \phi(v ) \in D_{n+1}$. If there exists $x$ such that $\neg \phi(\num{x}) \in D_k \cap T_k$, then $\neg \phi(\num{x}) \in D_n \cap T_n$, and so $\forall v \phi(v) \in D_{n+1}$.

	Now, we can prove the second part of our claim, namely that $D_{k+1} \cap T_{k+1} = D_{k+1} \cap T_{n+1}$. Again, we proceed by case distinction.
	
	\paragraph*{Atomic case -- arithmetical formulae}
	
	By definition if $\phi = (s=t)$ for $s,t \in \ClTerm_{\LPA}$, $\phi \in T_{n+1}$ iff $\phi \in T_{k+1}$ and $\phi$ is both in  $D_{k+1}$ and  $D_{n+1}$. 
	
	\paragraph*{Atomic case -- determinacy}
	
	Suppose that $\phi = D \psi$, with $\phi \in D_{k+1}$. Then $\psi \in D_k$ and by induction hypothesis, $\psi \in D_n$. In particular, by the atomic clauses of $\CT^-(D,T,T')$, we have both $\phi \in T_{k+1}$ and $\phi \in T_{n+1}$. 
	
	\paragraph*{Atomic case -- truth}
	
	Suppose that $\phi = T \psi$ and that $\phi \in D_{k+1} \cap T_{k+1}$. By the definition of $\mathscr{D}$, $\phi \in D_{k+1}$ iff $\psi \in D_k$ and by the definition of $\CT^-(D_k,T_k,T_{k+1})$, $\phi \in T_{k+1}$ iff $\psi in T_k$. So $\phi \in D_{k+1} \cap T_{k+1}$ iff $\psi \in D_k \cap T_k$ which by induction hypothesis is equivalent to $\psi \in D_k \cap T_n$, which again is equivalent to $\phi \in D_{k+1} \cap T_{n+1}$.

	\paragraph*{Negation}
	
	Suppose that $\phi = \neg \psi \in D_{k+1} \cap T_{k+1}$. Then $\psi \in D_k$ and $\psi \notin T_k$. By the induction hypothesis this is equivalent to $\psi \in D_n \setminus T_n$. The latter is equivalent to $\neg \psi \in D_{n+1} \cap T_{n+1}$.  
	
	\paragraph*{Conjunction}
	
	Suppose that $\phi = \psi \wedge \eta$. First assume that $\phi \in D_{k+1} \cap T_{k+1}$. We will show that $\phi \in T_{n+1}$. 
	
	If $\phi \in T_{k+1}$, then by the compositional conditions $\psi, \eta \in T_{k+1}$. Since $\phi \in D_{k+1}$, then either $\psi, \eta \in D_k$ or the negation of one of the conjuncts, say $\neg \psi$, is in $D_k \cap T_k$. The latter possibility cannot happen, since then by induction hypothesis we would have $\neg \psi \in T_{k+1}$ which violates the compositional conditions. 
	
	Since $\psi, \eta \in D_k$, we have $\psi, \eta \in T_{k}$, since otherwise by induction hypothesis we could not have $\psi, \eta \in T_{k+1}$. By induction hypothesis, $\psi, \eta \in T_{n+1}$, hence $\psi \wedge \eta \in T_{n+1}$. 
	
	Now assume that $\psi \wedge \eta \in D_{k+1}$ and $\psi \wedge \eta \notin T_{k+1}$. Then one of the conjuncts is determined and not true. Without loss of generality $\neg \psi \in D_k$ and $\neg \psi \in T_k$. Then by the induction hypothesis $\neg \psi \in T_{n+1}$ and by the compositional conditions $\psi \wedge \eta \notin T_{n+1}$. 
	
	\paragraph*{Universal quantifier}
	
	This is very similar to the previous case. Suppose first that $\forall x \phi(x) \in D_{k+1} \cap T_{k+1}$. Then for all $x$, $\phi(\num{x}) \in D_k$ and $\phi(\num{x}) \in T_{k}$, since otherwise for some $x$, $\neg \phi(\num{x}) \in D_k \cap T_k = D_k \cap T_{k+1}$, contradicting the compositional axioms. Hence, by the induction hypothesis, $\phi(\num{x}) \in T_{n+1}$ for all $x$, and consequently $\forall x \phi(x) \in T_{n+1}$. 
	
	Now, suppose that $\forall x \phi(x) \in D_{k+1} \setminus T_{k+1}$. In this case, there exists $x$ such that $\neg \phi(\num{x}) \in D_k$ and $\neg \phi(\num{x}) \in T_{k}$. By the induction hypothesis, $\neg \phi(\num{x}) \in D_k \cap T_{n+1}$, and so $\forall x \phi(x) \notin T_{n+1}$. 
	
This concludes the proof.
\end{proof}

\begin{claim} \label{cl_fixpoint}
The set $D_{\omega}$ is a fixpoint of the operator $\mathscr{D}(\cdot,T_{\omega})$, i.e.
\begin{displaymath}
	\mathscr{D}(D_{\omega},T_{\omega}) = D_{\omega}.
\end{displaymath}

\end{claim}
\begin{proof}
	We will verify by case distinction depending on the syntactic structure of a formula $\phi$ that if $\phi \in \mathscr{D}(D_{\omega}, T_{\omega})$, then $\phi \in  D_{\omega}$. The proof is entirely straightforward for any $\phi$ which is not a quantified statement. For the quantifier case, we will use the recursive saturation of our model.
	
	\paragraph*{Atomic case -- arithmetical formulae} Obvious.
	
	\paragraph*{Atomic case -- determinacy} Assume that $\phi = D(\psi)$ and that $\phi \in \mathscr{D}(D_{\omega},T_{\omega})$. This means that $\psi \in D_{\omega}$, so there exists a specific $k \in \omega$ such that $\psi \in D_k$, and thus $\phi \in D_{k+1}$, so $\phi \in D_{\omega}$. 
	
	\paragraph*{Atomic case -- truth} Suppose that $\phi = T(\psi)$. If $\phi \in \mathscr{D}(D_{\omega},T_{\omega})$, then $\psi \in D_k$ for some $k \in \omega$, and $\phi \in D_{k+1} \subset D_{\omega}$.
	
	\paragraph*{Negation} Obvious.
	
	\paragraph*{Conjunction} Suppose that $\phi: = \psi \wedge \eta \in \mathscr{D}(D_{\omega},T_{\omega}).$ Either $\psi, \eta \in D_{\omega}$ or (without loss of generality) $\neg \psi \in D_{\omega} \cap T_{\omega}$. 
	
	In the former case, we proceed straightforwardly: if $\psi \in D_k, \eta \in D_l$, then their conjunction is in $D_{m +1}$, where $m = \max(k,l)$.
	
	If $\neg \psi \in D_{\omega} \cap T_{\omega}$, then there exist $k,l$ such that $\neg \psi \in D_k$ and $\neg \psi \in D_l \cap T_l$ (and we can actually ignore the first part in the context of the proof). Then by the definition of the operator $\mathscr{D}$, $\psi \wedge \eta \in D_{l+1}$. 
	
	\paragraph*{Universal quantifier} Assume that $\phi = \forall v \psi(v)$ and that $\phi \in \mathscr{D}(D_{\omega}, T_{\omega})$. If there exists $x$ such that $\neg \psi(\num{x}) \in D_{\omega} \cap T_{\omega}$, then as in the previous case $\phi \in D_{\omega}$. 
	
	Now assume that $\psi(\num{x}) \in D_{\omega}$ for all $x$. We claim that there exists $k \in \omega $ such that for all $x$, $\psi(\num{x}) \in D_k$. Suppose otherwise. This means that the type:
	\begin{displaymath}
		\psi(\num{x}) \notin D_k, k \in \omega
	\end{displaymath}
	is finitely realised in $(M,I,T^*,D^*)$. Notice that the predicates $D_i$ are formally a single predicate $D_{v}(w)$, so the listed set of formulae actually only uses a finite part of the language. By recursive saturation, there exists an $x \in M$ such that 
	\begin{displaymath}
		\psi(\num{x}) \notin D_k, k \in \omega.
	\end{displaymath}
This contradicts our assumption. We conclude that for some $k$, $\psi(\num{x}) \in D_k$ and $\forall v \psi(v) \in D_{k+1}$. 
\end{proof}
	
\begin{claim} \label{cl_compositional}
	The predicate $T_{\omega}$ satisfies compositional conditions for the sentences in $D_{\omega}$. I.e., the following clauses are satisfied:
	\begin{itemize}
		\item $\forall s,t \in \ClTerm_{\LPA} \ \Big(T_{\omega}(s=t) \equiv \val{s} = \val{t} \Big).$
		\item $\forall x \Big(D_{\omega}(D(\num{x})) \rightarrow T_{\omega}D(\num{x}) \equiv D_{\omega}(x)\Big).$
		\item $\forall x \Big(D_{\omega}(T(\num{x})) \rightarrow T_{\omega}T(\num{x}) \equiv T_{\omega}(x)\Big).$
		\item $\forall \phi \in \Sent \ \Big( D_{\omega}(\neg \phi) \rightarrow T_{\omega}\neg \phi \equiv \neg T_{\omega} \phi\Big)$. 
		\item $\forall \phi \in \Sent \ \Big(D_{\omega}(\phi \wedge \psi) \rightarrow T_{\omega}(\phi \wedge \psi) \equiv \left(T_{\omega} \phi \wedge T_{\omega} \psi\right) \Big).$
		\item $\forall \phi, v \Big(\forall v \phi \in \Sent \cap D_{\omega} \rightarrow T_\omega \forall v \phi \equiv \forall x T_{\omega}\phi(\num{x})\Big).$ 
		\item $\forall \bar{s}, \bar{t} \in \ClTermSeq \forall \phi \ \Big(\bar{\val{s}} = \bar{\val{t}} \rightarrow T_{\omega}\phi(\bar{s}) \equiv T_{\omega} \phi(\bar{\val{t}}) \Big).$
	\end{itemize}
\end{claim}

\begin{proof}[Proof of Claim \ref{cl_compositional}]
	Since by definition $T_{\omega} = \bigcup_{n \in \omega} D_n \cap T_n$ and $D_{\omega} = \bigcup_{n \in \omega} D_n$. By induction on $n$, we will show that if $\phi \in D_{n}$, then compositional clauses hold for $\phi$. 
	
	If $n=0$, then the claim is trivial. So suppose now that the claim holds for all $k< n$. Let us check that compositionality holds, by case distinction, depending on the syntactic form of $\phi$.

	\paragraph*{Arithmetical atomic case.} If $\phi = (s=t)$ for some $s,t \in \ClTerm_{\LPA}$, then $\phi \in D_n$ for all $n \in \omega$ by definition and for all $n$, $T_n(s=t) \equiv \val{s} = \val{t}$ by the compositional clauses for the predicates $T_n$.

	\paragraph*{Atomic case -- determinacy}
	Let $\phi$ be a sentence of the form $D(\num{x})$ for some $x \in M$. We want to show that 
	\begin{displaymath}
		T_{\omega} D(\num{x}) \equiv D_{\omega}(x). 
	\end{displaymath}
	If $D_{\omega}(x)$ holds, say $x \in D_n$, then $D(\num{x}) \in D_{n+1}$ by the definition of $\mathscr{D}$, and $D(\num{x}) \in T_{n+1}$ by the definition of $\CT^-(D_n,T_n,T_{n+1}).$ 
	
	On the other hand, if $T_{\omega}(D(\num{x}))$ holds, then for some $n$, $D(\num{x}) \in T_n \cap D_n$. In particular (by the definition of $\CT^-(D_{n-1},T_{n-1},T_n)$), $x \in D_{n-1}$. Then $x \in D_{\omega}$ by definition.
	
	\paragraph*{Atomic case -- truth} Fix a sentence of the form $T(\num{x})$ for some $x \in M$ such that $T(\num{x}) \in D_{\omega}$. We want to check that the following equivalence holds:
	\begin{displaymath}
		T_{\omega}T(\num{x}) \equiv T_{\omega}(x).
	\end{displaymath}
	If $T(\num{x}) \in T_{\omega}$, there exists $n \in \omega$ such that $T(\num{x}) \in D_n \cap T_n$. Therefore, by the definition of $\mathscr{D}$, $x \in D_{n-1}$ and by the definition of $\CT^-(D_{n-1},T_{n-1},T_n)$, $x \in T_{n-1}$. Hence $x \in T_{\omega}$. 
	
	If, on the other hand, $x \in T_{\omega}$, then by definition, $x \in D_n \cap T_n$ for some $n \in \omega$. Then directly by the definitions of $\mathscr{D}$ and $\CT^-(D_n,T_n,T_{n+1})$, $T(\num{x}) \in T_{n+1} \cap D_{n+1}$. 

		\paragraph*{Negation} This step also does not depend on the induction hypothesis. We want to check the equivalence:
	\begin{displaymath}
		T_{\omega} \neg \phi \equiv \neg T_{\omega} \phi.
	\end{displaymath}
	Suppose that $T_{\omega} \neg \phi$ holds. Then there exists $k$ such that $T_{k} \neg  \phi$ and $D_k \neg \phi$ hold. By compositional clauses for $T_k$, we have:
	\begin{displaymath}
		T_k \neg \phi \equiv \neg T_k \phi.
	\end{displaymath}
	Suppose now that there exists $l$ such that $T_l \phi$ and $D_l \phi$ holds. By Claim \ref{cl_monotonicity}, if $l \leq k$, then $T_k \phi$ holds, and if $l \geq k$, then $T_l \neg \phi$ holds, in both cases contradicting compositionality. Therefore, there is none such $l$, and $\phi \notin T_{\omega}$. On the other hand, if $\phi \in D_n$ for some $n$ and $\phi \notin T_k \cap D_k$ for any $k < \omega$, then in particular $\phi \notin T_{n+1}$, and by the compositional clauses $\neg \phi \in D_{n+1} \cap T_{n+1}$.

	\paragraph*{Conjunction}
	Fix any sentence of the form $\phi \wedge \psi$. We want to show the following equivalence:
	\begin{displaymath}
		T_{\omega}(\phi \wedge \psi) \equiv T_{\omega} \phi \wedge T_{\omega} \psi.
	\end{displaymath}
	
	If $T_{\omega} \phi$, $T_{\omega} \psi$ hold, then there exist $n,k$ such that $\phi \in T_n \cap D_n$ and $\psi \in T_k \cap D_k$. By Claim \ref{cl_monotonicity}, assuming without loss of generality that $n \leq k$, we have $\phi, \psi \in T_k \cap D_k$, and, again by monotonicity, $\phi, \psi \in T_{k+1}$ and $\phi \wedge \psi \in D_{k+1}$ by definition. Then $\phi \wedge \psi \in T_{k+1} \cap D_{k+1}$ by compositionality. 
	
	Now suppose that $T_{\omega}(\phi \wedge \psi)$ holds. Then by definition, for some $k$, we have $T_{k}(\phi \wedge \psi)$ and $D_k(\phi \wedge \psi)$. By compositionality, $\phi \in T_k$ and $\psi \in T_k$. We want to check that $\phi \in D_k$ and $\psi \in D_k$, which would conclude the argument.
	
	Let $n+1$ be the earliest such that $\phi \wedge \psi \in D_{n+1}$. By definition of the operator $\mathscr{D}$, we either have $\phi, \psi \in D_n$ or (without loss of generality) $\neg \phi \in D_n, \neg \phi \in T_n$. If the former case holds, then $\phi, \psi \in D_k$. Suppose that the latter case holds. Then by Claim \ref{cl_monotonicity}, $\neg \phi \in T_k$ and by compositonality $\phi \notin T_k$ and so, again by compositionality, $\phi \wedge \psi \notin T_k$, contrary to our assumption.

	\paragraph*{Universal quantifier} 
	
	Fix a sentence $\forall v \phi(v) \in D_{\omega}$. We want to show that the following equivalence holds:
	\begin{displaymath}
		T_{\omega} \forall v \phi(v) \equiv \forall x \ T_{\omega}\phi(\num{x}).
	\end{displaymath}
	
	First, suppose that for all $x$, $\phi(\num{x}) \in D_{\omega} \cap T_{\omega}$. Then, by Claim \ref{cl_fixpoint}, $\forall v \phi(v) \in D_{\omega}$, hence there exists $n \in \omega$ such that $\forall  \phi(v) \in D_n$. We claim that $\forall (v) \phi(v) \in T_n$. Indeed, if this were not the case, then by compositionality we would have $\neg \phi(\num{x}) \in T_n$ for some $n$, and by Claim \ref{cl_monotonicity}, $\neg \phi(\num{x}) \in T_k$ for all $k \geq n$ which contradicts our assumption.
	
	Now assume that $\forall v \phi(v) \in D_{\omega} \cap T_{\omega}$. Then there exists $n$ such that $\forall v \phi(v) \in D_n \cap T_n$. Again, by compositionality it is enough to show that for all $x \in M$, $\phi (\num{x}) \in D_n$. Let $k+1$ be the earliest step such that $\forall v \phi(v) \in D_{k+1}$. Then either $\phi(\num{x}) \in D_k$ for all $x \in M$ (in which case we are done) or there exists $x$ such that $\neg \phi(\num{x}) \in D_k \cap T_k$. In this case $\neg \phi(\num{x}) \in D_n \cap T_n$ by monotonicity, and by compositionality $\phi(\num{x}) \notin T_n$ and hence, $\forall v \phi(v) \notin T_n$, contrary to the assumption.
	
	\end{proof}

Now, we can turn to the crucial part of the proof. Let $(M,D_{\omega},T_{\omega})$ be the model constructed before.

\begin{claim} \label{claim_extending_T_omega_to_T}
	The following theory is consistent:
	\begin{itemize}
		\item $\ElDiag(M,D_{\omega},T_{\omega})$.
		\item $\CT^-(M,D_{\omega},T_{\omega},T)$. (Compositionality).
		\item $\forall \phi \ \Big( D_{\omega}(\phi) \rightarrow T_{\omega}(\phi) \equiv T(\phi) \Big).$ (Compatibility).
	\end{itemize}
	
\end{claim}

The proof of this fact uses Enayat--Visser methods using Compositionality Claim \ref{cl_compositional}. Essentially, we want to find a fully compositional predicate $T$ extending a partial truth predicate $T_{\omega}$ which satisfies compositional clauses on a set of sentences $D_{\omega}$ satisfying strong enough closure conditions.  We will provide the proof of the claim in Appendix A. Now, fix $(M,D^{\bullet},T^{\bullet},T)$, a recursively saturated model of the above theory. Note that:
\begin{displaymath}
	(M,D^{\bullet},T^{\bullet}) \equiv (M,D_{\omega},T_{\omega}).
\end{displaymath}
We claim that $(M,D^{\bullet},T)$ is the desired model satisfying $\CD^-$. Let us verify that it indeed satisfies the axioms of this theory.

The clauses $\synt{T_1}, \synt{T_2}, \synt{T_4}, \synt{T_5}, \synt{T_6},$ and $\synt{R_1}$ are automatically satisfied in $(M,D^{\bullet},T)$, since $T$ satisfies $\CT^-$ over the arithmetical language extended with symbols $T,D$.

Let us check that $\synt{T_3}$ is satisfied. By the definition of $\CT(M,D^{\bullet},T^{\bullet},T)$ we have for all $\phi \in \Sent(M)$:

\begin{displaymath}
	(M,D^{\bullet},T^{\bullet},T) \models TT\phi \equiv T^{\bullet}\phi.
\end{displaymath}

Therefore, we need to check:
\begin{displaymath}
	(M,D^{\bullet},T^{\bullet,}T) \models D^{\bullet}(\phi) \rightarrow T^{\bullet}(\phi) \equiv T(\phi).
\end{displaymath} 
However, this follows immediately by the compatibility clause in the construction of $T^{\bullet},D^{\bullet},T$. Let us now check that the conditions for the predicate $D$ hold in $(M,D^{\bullet},T)$.

The conditions $\synt{D_1}, \synt{D_2}, \synt{D_3}, \synt{D_4}$ hold, by Fixed-point claim \ref{cl_fixpoint}, since $(M,D^{\bullet},T^{\bullet})$ is elementarily equivalent to $(M,D_{\omega},T_{\omega})$.

Let us check that $\synt{D_5}$ holds. Fix any $\phi,\psi \in M$. By Claim \ref{cl_fixpoint} and the elementary equivalence of $(M,D^{\bullet},T^{\bullet})$ with $(M,D_{\omega},T_{\omega})$, we have:
\begin{displaymath}
	D^{\bullet}(\phi \wedge \psi) \equiv \big((D^{\bullet}(\phi) \wedge D^{\bullet}(\psi)) \vee (D^{\bullet}(\neg \phi) \wedge T^{\bullet}(\neg \phi) ) \vee (D^{\bullet}(\neg \psi) \wedge T^{\bullet}(\neg \psi)) \big).
\end{displaymath} 
However, by $\synt{T_3}$, this is equivalent to:
\begin{displaymath}
	D^{\bullet}(\phi \wedge \psi) \equiv \big((D^{\bullet}(\phi) \wedge D^{\bullet}(\psi)) \vee (D^{\bullet}(\neg \phi) \wedge T(\neg \phi) ) \vee (D^{\bullet}(\neg \psi) \wedge T(\neg \psi)) \big),
\end{displaymath}
which proves that $\synt{D_5}$ holds.

Finally, let us check that $\synt{D_6}$ holds. Again, by Claim \ref{cl_fixpoint}:

\begin{displaymath}
	D^{\bullet}(\forall x \phi(x)) \equiv \Big(\forall x D^{\bullet}\phi(x) \vee \exists x \ D^{\bullet} \neg \phi(\num{x}) \wedge T^{\bullet} \neg \phi(\num{x}) \Big).
\end{displaymath} 
Like in the case of the axiom $\synt{D_5}$, by $\synt{T_3}$, this is equivalent to:
\begin{displaymath}
	D^{\bullet}(\forall x \phi(x)) \equiv \Big(\forall x D^{\bullet}\phi(x) \vee \exists x \ D^{\bullet} \neg \phi(\num{x}) \wedge T \neg \phi(\num{x}) \Big).
\end{displaymath} 	
This concludes the versification that $(M,D^{\bullet},T)$ satisfies the axioms of $\CD^-$ and, therefore, it concludes the whole proof.
\end{proof}

\section{Appendix A -- adapting the Enayat--Visser argument}

In proof of Theorem \ref{th_cdminus_conservative}, we have claimed that consistently, to the model $(M,D_{\omega},T_{\omega})$, a compositional truth predicate can be added which is compatible with the predicate $T_{\omega}$ on the sentences from $D_{\omega}$. This claim actually follows almost immediately from the original work of Enayat and Visser who showed that a model with a compositional satisfaction class defined on a set closed under direct subformulae can be elementarily extended and then expanded with a full satisfaction class containing the original one. 

Unfortunately, the class $T_{\omega}$ does not quite satisfy the above assumptions, since the set $D_{\omega}$ is not closed under direct subformulae. For instance, the following sentence is in $D_{\omega}$:
\begin{displaymath}
\exists x \ x = \qcr{0=0} \wedge T(x),
\end{displaymath}
whereas its subformula, $T(x)$, has instances which are not in $D_{\omega}$. Luckily, a natural weaker condition can be formulated which holds for $D_{\omega}$ and which guarantees that we can extend it to a full truth class in an elementary extension of a model.

For  technical convenience, the discussion below will be formulated in the language of satisfaction classes rather than truth predicates, as discussed in \cite{wcislo_definability_automorphisms}. We are inclined to say that the definition used there is standard. However, the actual definition used by various authors differ by slight details, so we refer to that work to be specific. In particular we will use the notion of \df{regular} satisfaction classes introduced in the Definition 13 of that article which we define below fo the convenience of the reader. 

In the formulation of this condition, we will make use of an admittedly technical notion of syntactic similarity. Both this notion and its variants were introduced and discussed in \cite{WcisloKossak}, \cite{lelyk_wcislo_local_collection}, and \cite{wcislo_definability_automorphisms} together with examples which might be helpful for the reader. 
\begin{definition} \label{def_syntactic similarity}
	Let $M \models \PA$ and let $\phi \in \form(M)$. Let $\alpha \in \Asn(\phi)$. By $\phi[\alpha]$ we mean the sentence resulting from $\phi$ by substituting a numeral $\num{\alpha(v)}$ for each variable $v$. 
	
	 Let $\phi,\psi \in \Sent(M)$ be two formulae. Let $\alpha \in \Asn(\phi), \beta \in \Asn(\psi)$ be two assignments. We say that the pairs $(\phi,\alpha)$ and $(\psi,\beta)$ are \df{syntactically similar} if there exists a formula $\eta$ and two sequences of terms $\bar{s},\bar{t} \in \ClTermSeq_{\LPA}$ such that:
	 \begin{itemize}
	 	\item $\bar{\val{s}} = \bar{\val{t}}.$ (The terms have equal values.)
	 	\item $\phi[\alpha] = \eta(\bar{s})$.
	 	\item  $\psi[\beta] = \eta(\bar{t})$.
	 \end{itemize}
 We denote this relation with $(\phi, \alpha) \sim (\psi,\beta)$. 
 
 We say that $\phi$ and $\psi$ are \df{syntactically similar} if for some $\alpha,\beta$ $(\phi,\alpha) \sim (\psi,\beta)$. We also denote this fact with $\phi \sim \psi$. 
\end{definition}

By a \df{full satisfaction class} on a model $M$ for a language $\mathscr{L}$, we mean a relation on $\form_{\mathscr{L}}(M) \times M$, such that for each $(\phi, \alpha) \in S$, $\alpha$ is a $\phi$-assignment (a function whose domain is the set of free variables of $\phi$) such that $S$ satisfies Tarski's compositional clauses for all $\phi \in \form_{\mathscr{L}}$. We say that $S$ is \df{regular} if $(\phi,\alpha) \sim (\psi,\beta)$ implies 
\begin{displaymath}
	S(\phi,\alpha) \equiv S(\psi,\beta).
\end{displaymath}
As we have mentioned, the details of the definition of a satisfaction class, and a discussion of some subtleties about that definition, can be found in \cite{wcislo_definability_automorphisms}. Now, we have to introduce a novel, somewhat technical notion.

\begin{definition} \label{def_locally_consistent}
Let $M \models \PA$. We say that a pair $D, S \subset M^2$ is \df{determinately compositional} if the following conditions hold:
\begin{itemize}
	\label{lem_ev_extendind_class_with_witnessing}
	\item[$\synt{D'1}$] $\forall s, t \in \Term_{\LPA} \forall \alpha \in \Asn (s=t) \ D(s=t,\alpha).$
	\item[$\synt{D'2}$] $\forall \phi \in \Sent \ D(T(\num{\phi}),\emptyset) \equiv D(\phi,\emptyset)$. 
	\item[$\synt{D'3}$] $\forall \phi \in \form, \alpha \in \Asn(\phi) \ D(D(\num{\phi}),\emptyset) \equiv D(\phi,\emptyset)$.
	\item[$\synt{D'4}$] $\forall \phi \in \form \forall \alpha \in \Asn(\phi)  D (\neg \phi, \alpha) \equiv D (
	\phi, \alpha). $
	\item[$\synt{D'5}$] $
	\forall \phi, \psi \in \form \forall \alpha \in \Asn(\phi \wedge \psi)\  D(\phi \vee \psi,\alpha) \equiv \\
	\left((D(\phi,\alpha) \wedge D(\psi,\alpha)) \vee (D(\phi,\alpha) \wedge S(  \phi,\alpha) \vee (D(\psi,\alpha) \wedge S ( \psi,\alpha) \right).$
	\item[$\synt{D'6}$] $			\forall \phi, v \ \Big(\phi \in \form, v \in \Var \forall \alpha \in \Asn(\phi) \ D(\exists v \phi) \equiv \\
	\left( \forall \beta \sim_v \alpha D(\phi,\beta) \vee \exists \beta \sim_v \alpha (D (\phi,\beta) \wedge S (\phi,\beta)) \right) \Big). $
	\item[$\synt{R'1}$] $\forall \phi,\psi \in \form \forall \alpha \in \Asn(\phi), \beta \in \Asn(\psi) \Big((\phi,\alpha) \sim (\psi,\beta) \rightarrow S(\phi,\alpha) \equiv S(\psi,\beta)\Big). $
	\item[$\synt{R'2}$] $\forall \phi,\psi \in \form \forall \alpha \in \Asn(\phi), \beta \in \Asn(\psi) \Big((\phi,\alpha) \sim (\psi,\beta) \rightarrow D(\phi,\alpha) \equiv D(\psi,\beta)\Big). $
	\item[$\synt{S1}$] $\forall s,t \in \Term_{\LPA} \forall \alpha \in \Asn(s,t)  \Big(S(s=t,\alpha) \equiv s^{\alpha} = t^{\alpha}\Big)$.
	\item[$\synt{S2}$] $\forall \phi \in \form, \alpha \in \Asn(\phi) \big(D(\phi,\emptyset)  \rightarrow S(D(\num{\phi}),\emptyset) \big).$
	\item[$\synt{S3}$] $\forall \phi \in \form, \alpha \in \Asn(\phi) \big(D(\phi,\emptyset) \rightarrow S(T(\num{\phi}),\emptyset) \rightarrow S(\phi,\emptyset)\big).$
	\item[$\synt{S4}$] $\forall \phi \in \form, \alpha \in \Asn(\phi) \big(D(\neg \phi,\alpha) \rightarrow S(\neg \phi,\alpha) \equiv \neg S(\phi,\alpha) \big).$
	\item[$\synt{S5}$] $\forall \phi, \psi \in \form, \alpha \in \Asn(\phi \vee \psi) \big(D(\phi \vee \psi, \alpha) \rightarrow S(\phi \vee \psi, \alpha) \equiv S(\phi,\alpha) \vee S(\psi,\alpha) \big).$
	\item[$\synt{S6}$] $\forall \phi \in \form, v \in \Var \big(D(\exists v \phi,\alpha) \rightarrow S(\exists v \phi, \alpha) \equiv \exists \beta \sim_v \alpha S(\phi,\beta)\big).$
\end{itemize}
\end{definition}

Whenever we write $\alpha \sim_v \beta$, we mean that $\alpha$ and $\beta$ agree as functions, except that we allow that $\alpha(v) \neq \beta(v)$ or that the value of $v$ is well defined only for one of the $\alpha$ or $\beta$. The expression $t^{\alpha}$ denotes the formally computed value of the terms $t$ under the assignment $\alpha$. Besides that, notice that in the above definition $\Sent$ and $\form$ denote formulae of the arithmetical language expanded with two \emph{unary} predicates $T(x), D(x)$. Even though in the proof we are using satisfaction classes for the technical convenience, the language under the satisfaction predicate remains unchanged. The axioms for determinate compositionality are very close to those of $\CD^-$. However, notice that we only require compositionality for the $S$ predicate whenever it applies to a determinate fact. Finally, for the technical convenience, we have moved from the language with conjunctions and the universal quantifiers to the one with disjunctions and the existential quantifiers. However, one can be straightforwardly translated into the other.

Observe that the axioms for determinate compositionality are also related to saying that a given class is compositional for a downward closed set of formulae. The crucial difference is that compositionality may obtain for some formulae \emph{applied only to some assignments}. We do not require that all atomic formulae satisfy disquotation and we do not require that the set of formulae on which $S$ is compositional is closed under taking subformulae. We might have an existential sentence which is determined which has some undetermined substitutional instances. However, detereminateness is closed enough that we have \emph{some} determined witness.

\begin{lemma} \label{lem_extending_determinate_compositionality}
Suppose that $(M_0,D_0,S_0)$ is a model of $\PA$ with a determinately compositional class. Then there exists an elementary extension $(M',D_0',S_0') \succeq (M,D_0,S_0)$ and $S' \supseteq S'$ such that $(M',S_0')$ is a full regular satisfaction class for the language $\LPA \cup \{D,S\}$.
\end{lemma}

Before we prove the above statement, let us make two observations.

\begin{cor}  \label{cor_T_extending_Tomega}
		Under the notation of the proof of Theorem \ref{th_cdminus_conservative}, Claim \ref{claim_extending_T_omega_to_T} follows.
\end{cor}
\begin{proof}
	Let $D_{\omega}$ and $T_{\omega}$ be the predicates constructed in the proof of Theorem 6. Let us define:
	\begin{eqnarray*}
		D_0(\phi,\alpha) & := & D_{\omega}(\phi^*[\alpha]) \\
		S_0(\phi,\alpha) & := & T_{\omega}(\phi^*[\alpha]),
	\end{eqnarray*}
	where by $\phi^*$ we mean $\phi$ with all conjunctions and universal quantifiers switched to disjunctions and the existential quantifiers via de Morgan laws. By Claims \ref{cl_fixpoint} and \ref{cl_compositional}, we obtain a determinately compositional pair. Now suppose that $S'$ is as in the conclusion of the Lemma. Then by defining:
	\begin{eqnarray*}
		T(\phi) & := & S'(\phi^{\circ},\emptyset), \\
		T^{\bullet}(\phi) & := & S_0'(\phi^{\circ},\emptyset) \\
		D^{\bullet}(\phi) & := & D_0'(\phi^{\circ},\emptyset),
	\end{eqnarray*} 
	(with $\cdot^{\circ}$  operation being the ``inverse'' of $\cdot^*$), we obtain a model $(M',T,T^{\bullet},D^{\bullet}) \models \CT^-(T,T^{\bullet},D^{\bullet})$ satisfying the Compatibility condition.
\end{proof}
\begin{remark} \label{rem_enayat_visser_bigger_languages}
	We can prove Fact \ref{fct_ctminus_w_wiekszych_jezykach} by the argument below, simply by removing the assumption that $(S_0,D_0)$ is a detereminately compositional pair in $M_0$ and by removing the Compatibility condition from the proof. 
\end{remark}

\begin{proof}[Proof of Lemma \ref{lem_extending_determinate_compositionality}]
In the proof, we will use a straightforward adaptation of the Enayat--Visser argument. Starting from $(M,D_0,S_0)$, we will construct a sequence of structures $(M_i,D_i,S^*_i,S_i)$ such that the chain $(M_i,D_i,S^*_i)$ is elementary, $S_{i+1}$ is compositional for the formulae in the model $M_i$, and $S_{i+1}$ agrees with $S_i$ on $M_{i-1}$ (setting $M_{-1} = \emptyset$ by convention) and with $S^*_i$ on $D_i$. 

So fix a model $(M_0,D_0,S^*_0)$ with a determinately compositional class. Let $S_0 = \emptyset$ and for each $i$, let $(M_{i+1},D_{i+1},S^*_{i+1},S_{i+1})$ be given as a model of the theory $\Theta_i$ with the following axioms:
\begin{itemize}
	\item $\ElDiag(M_i,D_i,S^*_i)$ (The elementary diagram).
	\item $\forall \phi,\psi \in \form \forall \alpha \in \Asn(\phi), \beta \in \Asn(\psi) \Big((\phi,\alpha) \sim (\psi,\beta) \rightarrow S_{i+1}(\phi,\alpha) \equiv S_{i+1}(\psi,\beta)\Big)$ (Regularity axiom).
	\item $\Comp(S_{i+1},\phi)$, where $\phi \in \form(M_i)$ (Compositionality scheme).
	\item $\forall \phi \in \form \forall \alpha \in \Asn(\phi) \Big(S^*_{i+1}(\phi,\alpha) \wedge D_{i+1}(\phi,\alpha) \rightarrow S_{i+1}(\phi,\alpha) \Big)$, where $\phi \in \form(M_i)$ (Compatibility axiom).
	\item $\forall \alpha \in \Asn(\phi) \Big(S_{i+1}(\phi,\alpha) \equiv S_i(\phi,\alpha)\Big)$, where $ \phi \in M_{i-1}$ (Preservation scheme). 
\end{itemize} 

By $\Comp(S_{i+1},\phi)$, we mean the finite set of axioms expressing that $S_{i+1}$ is compositional with respect to formula $\phi$, more precisely it is the disjunction of the following clauses:
\begin{itemize}
	\item $\exists s,t \in \Term_{\LPA} \ \Big(\phi = (s=t) \wedge \forall \alpha \in \Asn(s=t) \ S_{i+1}(s=t,\alpha) \equiv s^{\alpha} = t^{\alpha}\Big).$
	\item $\exists \psi \in \Sent \ \Big(\phi = D(\num{\psi}) \wedge S_{i+1}(\phi,\emptyset) \equiv D_{i+1}(\psi,\emptyset)\Big)$.
	\item $\exists \psi \in \Sent \ \Big(\phi = T(\num{\psi}) \wedge S_{i+1}(\phi,\emptyset) \equiv S^*_{i+1}(\psi,\emptyset)\Big)$.
	\item $\exists \psi \in \form \big( \phi = \neg \psi \wedge S_{i+1}(\phi,\alpha) \equiv \neg S_{i+1}(\psi,\alpha)\big).$
	\item $\exists \psi, \eta \in \form \big(\phi = (\psi \vee \eta) \wedge S_{i+1}(\phi,\alpha) \equiv S_{i+1}(\psi,\alpha) \vee S_{i+1}(\eta,\alpha)\big).$
	\item $\exists \psi \in \form, v \in \Var \big(\phi = \exists v \psi \wedge S_{i+1}(\phi,\alpha) \equiv \exists \beta \sim_v \alpha \ S_{i+1}(\psi,\beta)\big).$
\end{itemize}
We want to show that the theory $\Theta_i$ is consistent. We argue by compactness. Fix a finite theory $\Gamma \subset \Theta_i$. It is enough to find in $(M_i,D_i,S_i^*)$ a subset $S$ which satisfies the regularity and compatibility axioms and all the instances of compositionality and preservation schemes which appear in $\Gamma$.

Let $\phi_1, \ldots, \phi_n$ be all the formulae which appear under the $S_{i+1}$ predicate in those schemes. Without loss of generality, for any $i$ if $\phi_i$ is a disjunction, then either both of its disjuncts or neither of them appears among $\phi_j$.

Consider the relation $\unlhd$ defined on $M / \sim$, which is the transitive closure of the relation ,,$[\phi] \unlhd [\psi]$ iff there exists $\phi_0 \sim \phi$ and $\psi_0 \sim \psi$ such that $\phi$ is a direct subformula of $\psi$.'' Order the classes with respect to this relation. 

We define the predicate $S$ by induction on the the ordering finite $\unlhd$. For $[\phi]$ of the minimal rank, we let $(\phi,\alpha) \in S$ iff one of the following holds:
\begin{itemize}
	\item There exists a formula $\psi \in [\phi] \cap M_{n-1}$ and $S_i(\phi,\alpha)$ holds.
	\item $(\phi,\alpha) \in S^*_i \cap D_i$.
	\item There exist terms $t,s$ such that $\phi = (t=s)$ and $t^{\alpha} = s^{\alpha}$. (This condition is redundant if $D,S$ form a determinately compositional pair, but it has to be included in the proof of Fact \ref{fct_ctminus_w_wiekszych_jezykach}).
	\item There exists a term $t$ such that $\phi = Dt$ and $(t^{\alpha},\emptyset) \in D_i$. (This condition is redundant for $i>0$).
	\item There exists a term $t$ such that $\phi = Tt$ and $(t^{\alpha},\emptyset) \in S^*_i$. (This condition is redundant if $i>0$).  
\end{itemize} 
In effect, if a class $[\phi]$ is minimal, $\phi$ is not atomic, and does not intersect $M_{n-1}$ or $S_i^*$, we declare $\phi$ to be false under all assignments. 

For classes $[\phi]$ of higher rank, we define $S$ via compositional clauses. More specifically, for any $\psi \in [\phi]$ we let $(\psi,\beta) \in S$ iff one of the following conditions is satisfied:
\begin{itemize}
	\item There exists $\eta$ such that $\psi = \neg \eta$ and $S(\eta,\beta)$ doesn't hold.
	\item There exist $\eta, \zeta$ such that $\psi = \eta \vee \zeta$ and $S(\eta,\beta) \vee S(\zeta,\beta)$ holds.
	\item There exist $\eta \in \form, v \in \Var$ such that $\phi = \exists v \eta$ and for some $\gamma \sim_v \beta$, $S(\eta,\gamma)$ holds. 
\end{itemize}	
We will show that $S$ defined inductively as above satisfies $\Gamma$. 

\paragraph*{Compositionality and disquotation} We have defined $S$ inductively using the compositional clauses. Hence it satisfies compositional clauses for the non-atomic formulae independently of the choice of $S$ on the base formulae. The definition on formulae whose classes had minimal rank guarantees that $S$ satisfies the disquotational clauses for the atomic formulae.

\paragraph*{Regularity} The regularity axiom is satisfied, since we have defined the predicate $S$ simultaneously on the whole class of formulae of the minimal rank in a way which overtly respects syntactic similarity. Extending $S$ via compositional clauses respects regularity. 

\paragraph*{Preservation} For the formulae of the minimal rank, the predicate $S$ is defined so that $S(\phi,\alpha) \equiv S_i(\phi,\alpha)$ whenever $\phi$ is in the  $\sim$-class of a formula from $M_{i-1}$. By induction hypothesis, $S_i$ satisfies compositional clauses for formulae in $M_{i-1}$ and $M_{i-1}$ is closed under syntactic operations and under direct subformulae. Therefore $S_i$ agrees with $S$ on $M_{i-1}$.

\paragraph*{Compatibility}  By induction on the rank of the class of $\phi$, we will show that if $(\phi,\alpha) \in D_i$, then $(\phi, \alpha) \in S_i$ iff $(\phi,\alpha) \in S_i^*$. For the formulae of the minimal rank, this follows directly by the definition and by compatibility between $S_i$ and $S_i^*$. So let us assume that $[\phi]$ is not minimal and that $(\phi,\alpha) \in D_n$. 

If $\phi = \neg \psi$, then by the axioms for $D_n$, $(\psi, \alpha) \in D_n$, so by the induction hypothesis, $(\psi,\alpha) \in S_n^*$ iff $(\psi,\alpha) \in S$, so by the compositional clauses, the same holds for $\phi$. 

Assume that $\phi = \psi \vee \eta$. If both $(\psi,\alpha), (\eta,\alpha) \in D_n$, then the conclusion follows easily. Now suppose without loss of generality that  $(\psi,\alpha) \in D_n$ and that $( \psi, \alpha) \in S^*_n$. Then by the compositional clauses and the induction hypothesis, $(\phi,\alpha) \in S$ and $(\phi,\alpha) \in S^*_n$. 

Now suppose that $\phi = \exists v \psi$. Again, if for all $\beta \sim_v \alpha$, we have $(\psi,\beta) \in D_n$, then the conclusion follows immediately by the compositional clauses, and if $(\psi,\beta) \in D_n \cap S^*_n$ for some $\beta$, then $(\phi, \alpha) \in S^*_n$ iff $(\phi,\alpha) \in S$.

This concludes the proof of the compatibility scheme, and of the consistency of $\Theta_i$. 

\paragraph*{The conclusion of the proof}

Finally, let us set $(\phi,\alpha) \in S$ iff there exists $n$ such that $\phi \in M_n$ and $(\phi,\alpha) \in S_{n+1}$. We check that by the preservation and compositionality schemes, the obtained predicate is compositional. Thanks to the compatibility scheme, it agrees with $S^* = \bigcup_{n \in \omega} S^*_n$ on $D' = \bigcup_{n \in \omega} D_n$. 
\end{proof}

\section{Appendix B -- Formalised syntactic notions} \label{sec_app_syntax}

Throughout the paper, we have made use of a number of syntactic notions. For the convenience of the reader, we gather them in one place.

\begin{itemize}
	\item $\alpha \sim_v \beta$ is a formula expressing ``$\alpha$ and $\beta$ are functions and they are equal, possibly except that $v \in \dom(\alpha) \setminus \dom(\beta)  \cup \dom(\beta) \setminus \dom(\alpha)$ or $\alpha(v) \neq \beta(v)$.''
	\item $\Asn(\phi)$. If $\phi$ is a formula, by $\Asn(\phi)$, we denote the set of assignments for $\phi$ (that is, functions, whose domain is the set of free variables of $\phi$).
	\item $\ClTerm_{\LPA}(x)$ is a formula expressing ``$x$ is a closed arithmetical term.''
	\item $\ClTermSeq_{\LPA}(x)$ is a formula expressing ``$x$ is a sequence of closed arithmetical terms.''
	\item $\form(x)$ is a formula expressing ``$x$ is a formula of the language $\LPA \cup \{D,T\}$.''
	\item $\Sent(x)$ is a formula expressing ``$x$ is a sentence of the language $\LPA \cup \{D,T\}$.''
	\item $\val{s}, \val{t}$. If $t$ is (a code of) a closed arithmetical term, by $\val{t}$, we denote the formally computed value of $t$. (As we have mentioned in Section \ref{sec_prelim}, we treat $\val{t}$ as a self-standing expression). If $\bar{t}$ is a sequence of closed arithmetical terms, then $\bar{\val{t}}$ denotes the sequence of their values. 
	\item $t^{\alpha}, s ^{\alpha}$. If $t$ is (a code of) an arithmetical term and $\alpha$ is an assignment, by $t^{\alpha}$ we denote the formally computed value of $t$ under this assignment.
	\item $\num{x}$. By $\num{x}$ we mean some canonically chosen closed term with the formally computed value $x$.  
\end{itemize}

\section*{Acknowledgments}

This research was supported by the NCN MAESTRO grant 2019/34/A/HS1/00399 ``Epistemic and Semantic Commitments of Foundational Theories.''

\end{document}